\documentclass{gtmon_a}
\pdfoutput=1

\usepackage{graphicx}
\usepackage[all]{xy}

\proceedingstitle{Groups, homotopy and configuration spaces (Tokyo
  2005)}
\conferencestart{5 July 2005}
\conferenceend{11 July 2005}
\conferencename{Groups, homotopy and configuration spaces, 
                in honour of Fred Cohen's 60th birthday}
\conferencelocation{University of Tokyo, Japan}

\editor{Norio Iwase}
\givenname{Norio}
\surname{Iwase}

\editor{Toshitake Kohno}
\givenname{Toshitake}
\surname{Kohno}

\editor{Ran Levi}
\givenname{Ran}
\surname{Levi}

\editor{Dai Tamaki}
\givenname{Dai}
\surname{Tamaki}

\editor{Jie Wu}
\givenname{Jie}
\surname{Wu}

\title{Computation of the homotopy of the spectrum \texttt{tmf}}
\author{Tilman Bauer}
\givenname{Tilman}
\surname{Bauer}
\address{Mathematisches Institut der Universit\"at M\"unster\\\newline
Einsteinstr. 62\\48149 M\"unster\\Germany}
\email{tbauer@uni-muenster.de}
\urladdr{http://wwwmath.uni-muenster.de/u/tbauer/}

\volumenumber{13}
\issuenumber{}
\publicationyear{2008}
\papernumber{02}
\startpage{11}
\endpage{40}

\doi{}
\MR{}
\Zbl{}

\arxivreference{math/0311328}


\received{30 September 2005}
\revised{1 September 2006}
\accepted{11 October 2006}
\published{22 February 2008}
\publishedonline{22 February 2008}
\proposed{}
\seconded{}
\corresponding{}
\version{}

\let\xysavmatrix\xymatrix
\def\xymatrix{\disablesubscriptcorrection\xysavmatrix}
\AtBeginDocument{\let\bar\wbar\let\tilde\wtilde\let\hat\what}

\newcommand {\h}{\mathfrak h}

\renewcommand {\S} {\mathbf S}

\newcommand {\F}{\mathbf F}

\renewcommand{\Z}{\mathbf Z}

\newcommand{\id}{\operatorname{id}}

\newcommand {\sm}{\wedge}

\swapnumbers

\makeatletter
\def\cnewtheorem#1[#2]#3{\newtheorem{#1}{#3}[section]
\expandafter\let\csname c@#1\endcsname\c@subsection}

\cnewtheorem{lemma}[subsection]{Lemma}
\cnewtheorem{thm}[subsection]{Theorem}
\cnewtheorem{corollary}[subsection]{Corollary}
\cnewtheorem{prop}[subsection]{Proposition}
\cnewtheorem{numb}[subsection]{} 
\theoremstyle{definition}
\cnewtheorem{example}[subsection]{Example}
\newtheorem*{defn}{Definition} 
\cnewtheorem{remark}[subsection]{Remark}
  \let\c@equation\c@subsection
  \def\theequation{\thesubsection}
\makeatother

\DeclareMathOperator{\Ext}{Ext}

\DeclareMathOperator{\Sq}{Sq}
\newcommand{\tmf}{\mathtt{tmf}}

\def\smashop#1_#2{%
\displaystyle{#1_{%
\hbox to 0pt{\hss$\scriptstyle{#2}$\hss}}\;}}
\makeatother
\newdir { >} {{}*!/-5pt/@{>}}

\newcommand{\ZZ}{\mathbf{Z}}
\renewcommand{\SS}{\mathbf{S}}

\newcommand{\M}{\mathcal{M}}

\def\fourmatricmassey#1#2#3#4{
\left\langle \begin{pmatrix} #1 & #2 \end{pmatrix},
\begin{pmatrix} #3 & #4 \\ #4 & #3 \end{pmatrix},
\begin{pmatrix} #1 & #2 \\ #2 & #1 \end{pmatrix},
\begin{pmatrix} #3 \\ #4 \end{pmatrix} \right\rangle
}
\newcommand{\undl}[1]{\underline{\smash{#1}}}
\newcommand{\oneover}[1]{\hbox{$\frac1{#1}$}}

\newsavebox{\kbar}
\sbox{\kbar}{$\Bar\kappa$}

\begin{document}

\begin{htmlabstract}
This paper contains a complete computation of the homotopy ring of the
spectrum of topological modular forms constructed by Hopkins and
Miller. The computation is done away from 6, and at the
(interesting) primes 2 and 3 separately, and in each of the latter
two cases, a sequence of algebraic Bockstein spectral sequences is
used to compute the E<sub>2</sub> term of the elliptic Adams&ndash;Novikov spectral
sequence from the elliptic curve Hopf algebroid. In a further step,
all the differentials in the latter spectral sequence are
determined. The result of this computation is originally due to
Hopkins and Mahowald (unpublished).
\end{htmlabstract}

\begin{abstract}
This paper contains a complete computation of the homotopy ring of the
spectrum of topological modular forms constructed by Hopkins and
Miller. The computation is done away from $6$, and at the
(interesting) primes $2$ and $3$ separately, and in each of the latter
two cases, a sequence of algebraic Bockstein spectral sequences is
used to compute the $E_2$ term of the elliptic Adams--Novikov spectral
sequence from the elliptic curve Hopf algebroid. In a further step,
all the differentials in the latter spectral sequence are
determined. The result of this computation is originally due to
Hopkins and Mahowald (unpublished).
\end{abstract}

\begin{asciiabstract}
This paper contains a complete computation of the homotopy ring of the
spectrum of topological modular forms constructed by Hopkins and
Miller. The computation is done away from 6, and at the (interesting)
primes 2 and 3 separately, and in each of the latter two cases, a
sequence of algebraic Bockstein spectral sequences is used to compute
the E_2 term of the elliptic Adams-Novikov spectral sequence from the
elliptic curve Hopf algebroid. In a further step, all the
differentials in the latter spectral sequence are determined. The
result of this computation is originally due to Hopkins and Mahowald
(unpublished).
\end{asciiabstract}

\maketitle

\section{Introduction}

In \cite{hopkins,hopkins-mahowald,hopkins-miller}, Hopkins, Mahowald,
and Miller constructed a new homology theory called $\tmf$, or
topological modular forms (it was previously also called $eo_2$). This
theory stands at the end of a development of several theories called
elliptic cohomology (Landweber, Ravenel and Stong
\cite{Landweber-Ravenel-Stong} Landweber
\cite{Landweber-elliptic-LNM}, Franke\cite{Franke:elliptic} and
others). The common wish was to produce a universal elliptic
cohomology theory, ie, a theory $\operatorname{Ell}$ together with an
elliptic curve $C_u$ over $\pi_0\operatorname{Ell}$ such that for any
complex oriented spectrum $E$ and any isomorphism of the formal group
associated to $E$ with the formal completion of a given elliptic curve
$C$, there are unique compatible maps $\operatorname{Ell} \to E$ and
$C_u \to C$. It was immediately realized that this is impossible due
to the nontriviality of the automorphisms of the ``universal''
elliptic curve. Earlier theories remedied this by inverting some small
primes (eg, by considering this moduli problem over $\ZZ[1/6]$), or
by considering elliptic curves with level structures, or both. The
Hopkins--Miller approach, however, is geared in particular towards
studying the small prime phenomena that arise; as a trade-off, their
theory $\tmf$ is not complex orientable. This, however, turned out to
have positive aspects. Hopkins and Mahowald first computed the
homotopy groups of $\tmf$ and realized that the Hurewicz map
$\pi_*\SS^0 \to \pi_*\tmf$ detects surprisingly many classes, at least
at the primes $2$ and $3$; thus $\tmf$ was found to be a rather good
approximation to the stable stems themselves at these primes. Their
computation was never published, and the aim of this note is to give a
complete calculation (in a different way than theirs) of the homotopy
of $\tmf$ at the primes $2$ and $3$.

Along with the construction of $\tmf$, which is constructed as an $A_\infty$--ring spectrum and is shown to be $E_\infty$, Hopkins and Miller set up an Adams--Novikov type spectral sequence that converges to $\pi_*(\tmf)$; its $E_2$ term is given as an Ext ring of a Hopf algebroid,
\[
E_2^{s,t} = \Ext^{s,t}_{(A,\Gamma)}(A,A),
\]
where all the structure of the Hopf algebroid $(A,\Gamma)$ can be
described completely explicitly. One particular feature of this Hopf
algebroid is that both $A$ and $\Gamma$ are polynomial rings over
$\ZZ$ on \emph{finitely} many generators; it is this finiteness that
generates a kind of periodicity in $\pi_*(\tmf)$ and makes it possible
to compute the whole ring and not only the homotopy groups up to a
certain dimension.

In \fullref{homalghopfalg} I will review some notions and tools from
the homological algebra of Hopf algebroids; in \fullref{ellalg}, the
elliptic curve Hopf algebroid is defined and studied. The short
\fullref{largeprimes} contains the computation of $\pi_*\tmf$ away
from $6$; Sections \ref{extthree} and \ref{differentialsatthree}
contain the computation of the Ext term and the differentials,
respectively, of the spectral sequence at the prime $3$; in Sections
\ref{exttwo} and \ref{differentialsattwo}, the same is done for the
prime $2$.

\subsubsection*{Acknowledgements}
I do not claim originality of any of the results of this paper. As
already noted, the computation was first done by Hopkins and
Mahowald. Rezk \cite{Rezk:math512} also has some lecture notes with an
outline of the computations that need to be done when computing the
homotopy from the Hopf algebroid. I found it worthwhile to make this
computation available anyway since it is, in my opinion, an
interesting computation, and explicit knowledge of it is useful when
using $\tmf$ as an approximation to stable homotopy theory.

I am very grateful to Doug Ravenel, Andy Baker, John Rognes, and the
anonymous referee for suggesting a number of improvements and
corrections of this computation. Most importantly, I am deeply
indepted to Mike Hopkins, my thesis advisor, from who I learned
everything I know today about $\tmf$.

The \LaTeX\ package {\tt sseq} which produced the numerous spectral
sequence charts in this paper, is available from CTAN servers or from
\href{http://wwwmath.uni-muenster.de/u/tbauer/}{the author's home
page}.

\section{Homological algebra of Hopf algebroids} \label{homalghopfalg}

In this section I review some important constructions that help in computing the cohomology of Hopf algebroids. All of this is well-known to the experts, but explicit references are a little hard to come by.

\begin{defn}
The category of \textit{Hopf algebroids} over a ring $k$ is the category of co\-groupoid objects $(A,\Gamma)$ in the category of commutative $k$--algebras such that the left (equivalently, right) unit map $\eta_L\co A \to \Gamma$ (resp.\ $\eta_R$) is flat.

The category of left \textit{comodules} over a Hopf algebroid $(A,\Gamma)$ is the category of $(A,\Gamma)$--left coaction objects in $k$--modules.

(For a more explicit description, consult Ravenel \cite[Appendix A1]{Ravenel:green2}.)
\end{defn}

By convention, we will always consider $\Gamma$ as an $A$--module using the left unit $\eta_L$ unless otherwise stated.
The flatness condition (which not all authors consider part of the definition of a Hopf algebroid) ensures that the category of comodules is an abelian category, which has enough injectives so that homological algebra is possible. We will abbreviate
\[
H^n(A,\Gamma;M) = \Ext^n_{(A,\Gamma)-\text{comod}}(A,M).
\]
and $H^*(A,\Gamma) = H^*(A,\Gamma; A)$. In our situation, all objects are $\Z$-graded, and we write
\[
H^{n,m}(A,\Gamma; M) = \Ext^n_{(A,\Gamma)-\text{comod}}(A,M[-m])
\]
where $(M[m])_n=M_{n-m}$.

\begin{defn}
A \textit{natural transformation} between two maps of Hopf algebroids 
\[
f=(f_0,f_1),\;g=(g_0,g_1)\co (A,\Gamma) \to (A',\Gamma')
\]
is a natural transformation of functors between cogroupoid objects; explicitly, it is an algebra map $c\co \Gamma \to A'$ such that $c\circ\eta_L = f_0$ and $c\circ\eta_R = g_0$, and such that the following two composites are equal:
\[
\xymatrix@R=0pt{
&& A' \otimes_A \Gamma \ar[dr]\\
\Gamma \ar[r]^(0.35){\psi} & \Gamma \otimes_A \Gamma \ar[ur]^-{c\otimes\id} \ar[dr]_-{\id \otimes c} & & \Gamma',\\
&& \Gamma \otimes_A A' \ar[ur]
}
\]
where $\psi$ is the Hopf algebroid diagonal. In the lower row, $A'$ is an $A$--module by means of the map $f_0$, and the rightmost map is induced by $f_*\co \Gamma \otimes_A A' \to \Gamma'\otimes_{A'} A'$. The upper row is defined analogously with $g$.

An \textit{equivalence} of Hopf algebroids $(A,\Gamma) \simeq (A',\Gamma')$ consists of maps 
\[
f\co (A,\Gamma) \leftrightarrows (A',\Gamma'):\!g
\]
and natural transformations between the identity and $f\circ g$ and between the identity and $g \circ f$. 
\end{defn}

\begin{remark}
The condition that $c \circ \eta_L = \id$ for a natural transformation from the identity to $f\circ g$ is saying that $c$ is an $A$--algebra map.
\end{remark}

An equivalence of Hopf algebroids induces an equivalence of comodule categories; in particular,
\[
H^*(A,\Gamma; M) \cong H^*(A',\Gamma'; f_*M).
\]

\begin{remark}
In the language of stacks, this has the following interpretation. To
every Hopf algebroid $(A,\Gamma)$ there is associated a stack
$\M_{(A,\Gamma)}$; this assignment is left adjoint to the inclusion of
stacks into groupoid valued functors. Equivalent Hopf algebroids give
rise to equivalent associated stacks. The Hopf algebroid cohomology
with coefficients in the comodule $M$ can be interpreted as the
cohomology of the quasi-coherent sheaf $\tilde M$ on the stack.
\end{remark}

\subsection{Base change}

Let $(A,\Gamma)$ be any Hopf algebroid, and let $A \to A'$ be a morphism of rings. Define $\Gamma' = A' \otimes_A \Gamma \otimes_A A'$, where the left tensor product is built using the map $\eta_L\co A \to \Gamma$, and the right one using $\eta_R\co A \to \Gamma$. Then $(A',\Gamma')$ is also a Hopf algebroid, and there is a canonical map of Hopf algebroids $(A,\Gamma) \to (A',\Gamma')$.

The following theorem is proved by Hovey and Sadofsky \cite[Theorem
  3.3]{Hovey-Sadofsky:invertible} and Hovey \cite[Corollary 5.6]{Hovey:morita}. The authors attribute it to Hopkins.

\begin{thm}\label{basechg}
If $(A,\Gamma) \to (A',\Gamma')$ is such a map of Hopf algebroids, and there exists a ring $R$ and a morphism $A' \otimes_A \Gamma \to R$ such that the composite
\[
A \xrightarrow{1 \otimes \eta_R} A' \otimes_A \Gamma \to R
\]
is faithfully flat, then it induces an equivalence of comodule categories, and in particular an isomorphism
\[
H^*(A,\Gamma) \xrightarrow{\cong} H^*(A',\Gamma')
\]
\end{thm}

This theorem is reminiscent of the change of rings theorem \cite[A1.3.12]{Ravenel:green2}. In fact, the change of rings theorem is enough to prove all the consequences of the above theorem used in this work, but I will use \fullref{basechg} for expository reasons. 

\subsection{The algebraic Bockstein spectral sequence}

Let $I \triangleleft (A,\Gamma)$ be an invariant ideal of a Hopf
algebroid, ie, $I \triangleleft A$ and $\eta_R(I) \subseteq I
\Gamma$. Denote by $(A,\Gamma)/I$ the induced Hopf algebroid $(A/I, I
\backslash \Gamma / I)$.

\begin{thm}[Miller, Novikov] \label{algbocksteinss}
There is a spectral sequence of algebras, called the algebraic Bockstein spectral sequence or the algebraic Novikov spectral sequence:
\[
E_2 = H^*\left(\left(A,\Gamma\right)/I;I^n/I^{n+1}\right) \Longrightarrow H^*(A,\Gamma)
\]
arising from the filtration of Hopf algebroids induced by the filtration
\[
0 < A/I < A/I^2 < \dots < A
\]
If $A$ is Cohen--Macaulay and $I \triangleleft A$ is a regular ideal, then $I^n/I^{n+1} \cong \operatorname{Sym}^n(I/I^2)$.
\end{thm}

In the applications in this paper, $I=(x)$ will always be a principal ideal, and $\operatorname{Sym}^*(I/I^2) = A'[y]$ will be a polynomial algebra. Thus the $E_2$ term simplifies to
\[
E_2 = H^*\left(\left(A,\Gamma\right)/I\right)[y].
\]
If $x=p$, then this is the ordinary (bigraded) Bockstein spectral sequence.
In the topological context, Hopf algebroids are usually graded, and thus the spectral sequence is really tri-graded. In displaying charts, I will disregard the homological filtration degree in the above spectral sequence.

\section{The elliptic curve Hopf algebroid} \label{ellalg}

The Hopf algebroid of Weierstrass elliptic curves and their (strict) isomorphisms is given by:
\begin{gather*}
A = \ZZ[a_1,a_2,a_3,a_4,a_6]; \quad |a_i| = 2i;\\
\Gamma = A[r,s,t]; \quad |r|=4;\; |s|=2;\; |t|=6
\end{gather*}
and the following structure maps:

The left units $\eta_L$ are given by the standard inclusion $A \hookrightarrow 
\Gamma$, the right units are:
\begin{align*}
\eta_R(a_1) &= a_1 + 2s \\
\eta_R(a_2) &= a_2 - a_1 s + 3 r - s^2\\
\eta_R(a_3) &= a_3 + a_1 r + 2 t\\
\eta_R(a_4) &= a_4 - a_3 s + 2 a_2 r - a_1 t - a_1 r s - 2 s t + 3 r^2\\
\eta_R(a_6) &= a_6 + a_4 r - a_3 t + a_2 r^2 - a_1 r t - t^2 + r^3.
\end{align*}
The comultiplication map is given by:
\begin{align*}
\psi(s) &= s \otimes 1 + 1 \otimes s\\
\psi(r) &= r \otimes 1 + 1 \otimes r\\
\psi(t) &= t \otimes 1 + 1 \otimes t + s \otimes r.
\end{align*}
This Hopf algebroid classifies plane cubic curves in Weierstrass form; the universal Weierstrass curve is defined over $A$ in affine coordinates $(x,y)$ as
\begin {equation}
y^2+a_1 x y + a_3 y = x^3 + a_2 x^2 + a_4 x + a_6, \label{weierstrass}
\end {equation}
and the universal strict isomorphism classified by $\Gamma$ is the coordinate change
\begin{align*}
x &\mapsto x+r\\
y &\mapsto y+sx+t.
\end{align*}
It is easy, and classical, to compute the ring of invariants of this Hopf algebroid:
\[
H^{0,*}(A,\Gamma) = \Z[c_4,c_6,\Delta]/(c_4^3-c_6^2-1728\Delta)
\]
where the classes $c_4$, $c_6$, $\Delta$ are given in the universal
case, and hence for every Weierstrass curve, by the equations (cf
Silverman \cite[III.1]{Silverman:arithmetic})
\begin{align*}
b_2 &= a_1^2+4a_2\\
b_4 &= 2a_4+a_1a_3\\
b_6 &= a_3^2 + 4a_6\\
c_4 &= b_2^2-24b_4\\
c_6 &= -b_2^3+36b_2b_4-216 b_6\\
\Delta &= \oneover{1728}(c_4^3-c_6^2).
\end{align*}
A cubic of the form \eqref{weierstrass} is an elliptic curve if and only if $\Delta$, the discriminant, is invertible. An elliptic curve is a $1$--dimensional abelian group scheme, and thus the completion at the unit (which is the point at infinity) is a formal group. In fact, this construction even gives a formal group when the discriminant is not invertible because the point at infinity in a Weierstrass curve is always a smooth point. Thus we get a map of Hopf algebroids
\[
H\co (MU_*,MU_*MU) \to (A,\Gamma).
\]
and a map from the $MU$--based Adams--Novikov spectral sequence to the
elliptic spectral sequence, which converges to the Hurewicz map
$h\co \SS^0 \to \tmf$. This last statement is due to the fact that
there is a complex oriented ring spectrum $\underline A$ with
$\pi_*\underline A = A$; it can be constructed as $\underline A = \tmf
\sm Y(4)$, where $Y(4)$ is the Thom spectrum of $\Omega U(4) \to \Z
\times BU$ (Rezk \cite[Sections 12--14]{Rezk:math512}).

Let $MU_*=\Z[x_1,x_2,\dots]$ with $|x_i|=2i$, and $MU_*MU = MU_*[m_1,m_2,\dots]$ with $|m_i|=2i$.
By direct verification, we may choose the generators $x_i,\; m_i$ such that
\begin{equation} \label{ANSScomparison}
\begin{array}{cc}
H(x_i) = a_i \quad \text{ for } i=1,2;\\
H(m_1) = s \quad \text{ and } \quad H(m_2) = r.
\end{array}
\end{equation}
This very low-dimensional observation enables us to compare the $\tmf$ spectral sequence with the usual $MU$--Adams--Novikov spectral sequence and will be the key to determining the differentials in sections \ref{differentialsatthree} and \ref{differentialsattwo}.

\section{The homotopy of \texorpdfstring{$\tmf$}{tmf} at primes \texorpdfstring{$p>3$}{p>3}} \label{largeprimes}

As was mentioned in the introduction, $\tmf$ becomes complex oriented at large primes, and thus the homotopy ring becomes very simple; it is isomorphic to the ring of classical modular forms, which is generated by the Eisenstein forms $E_4$ and $E_6$.

We first state a result valid whenever $2$ is invertible in the ground ring.

Let $\tilde A = \Z[\oneover 2,a_2,a_4,a_6]$, $f\co A[\oneover 2] \to \tilde A$ be the obvious projection map sending $a_1$ and $a_3$ to $0$, and
\[
\tilde \Gamma = \tilde A \otimes_A \Gamma \otimes_A \tilde A = A[\oneover 2,r,s,t]/(a_1,a_3,\eta_R(a_1),\eta_R(a_3)) = \tilde A[r].
\]

\begin{lemma} \label{equivalenceawayfromtwo}
The Hopf algebroids $(A[\oneover 2],\Gamma[\oneover 2])$ and $(\tilde A, \tilde \Gamma)$ are equivalent.
\end{lemma}
\begin{proof}
Define a map $g\co (\tilde A, \tilde \Gamma) \to (A[\frac12],\Gamma[\frac12])$ in the opposite direction of $f$ by \begin{align*}
g(a_2) &= a_2+\oneover 4 a_1^2,\\
g(a_4) &= a_4 + \oneover 2 a_1 a_3,\\
g(a_6) &= a_6 + \oneover 4 a_3^2, \text{ and}\\
g(r) &=r.
\end{align*}
Then clearly $f \circ g = \id_{(\tilde A,\tilde \Gamma)}$, and if we define $c\co \Gamma[\oneover 2] \to A[\oneover 2]$ to be $A$--linear and 
\[
c(r) = 0; \quad c(s) = -\oneover 2 a_1; \quad c(t) = -\oneover 2 a_3,
\]
we find that $c \circ \eta_R = g \circ f$.
\end{proof}

\begin{remark} \label{completethecube}
This proof is the formal version of the following algebraic argument: As $2$ becomes invertible, we can complete the squares in the Weierstrass equation \eqref{weierstrass} by the substitution
\[
y \mapsto y-\oneover 2 a_1 x - \oneover 2 a_3.
\]
That is, we can find new parameters $x$ and $y$ such that the 
elliptic curve is given by
\begin{equation} \label{hyperellipticeq}
y^2 = x^3 + a_2 x^2 + a_4 x + a_6.
\end{equation}
All the strict isomorphisms of such curves are given by substitutions $x \mapsto x+r$, hence the Hopf algebroid is equivalent to $(\tilde A, \tilde \Gamma)$.
\end{remark}

When $3$ is also invertible in the ground ring, we can in a similar fashion complete the cube on the right hand side of \eqref{hyperellipticeq} by sending $x \mapsto x-\frac13 a_2$ and obtain
\begin{prop}\

\begin{enumerate}
\item \label{largeprimeequiv} The Hopf algebroid $\left(A[\frac16],\Gamma[\frac16]\right)$ is equivalent to the discrete Hopf algebroid $(A',\Gamma'=A')$, where $A' = \Z[\frac16,a_4,a_6]$.
\item $\displaystyle{
H^{n,*}(A[\textstyle{\frac16}],\Gamma[\textstyle{\frac16}]) = \begin{cases} \Z[\frac16,a_4,a_6]; & n=0\\0; & \text{otherwise} \end{cases}
}$ \label{largeprimecohom}
\item \label{largeprimehomotopy} $\pi_*\left(\tmf[\frac16]\right) = \Z[\frac16,a_4,a_6]$.
\end{enumerate}
\end{prop}
\begin{proof}
Clearly, \eqref{largeprimeequiv} implies \eqref{largeprimecohom}, and the latter implies \eqref{largeprimehomotopy} because the elliptic spectral sequence collapses with $E_\infty$ concentrated on the horizontal axis.
\end{proof}

\section{The cohomology of the elliptic curve Hopf algebroid at\break \texorpdfstring{$p=3$}{p=3}} \label{extthree}

As a warm-up for the more difficult case of $p=2$, we first compute $H^{**}(A_{(3)},\Gamma_{(3)})$.

\fullref{equivalenceawayfromtwo} gives us an equivalence of $(A_{(3)},\Gamma_{(3)})$ with $(\tilde A, \tilde \Gamma)$, where
\[
\tilde A = \Z_{(3)}[a_2,a_4,a_6] \quad \text{and} \quad \tilde \Gamma = \tilde A[r].
\]
To simplify the Hopf algebroid further, we map $\tilde A$ to $A' = \Z_{(3)}[a_2,a_4]$ by sending $a_6$ to $0$. Now
\[
g\co \tilde A \xrightarrow{1 \otimes \eta_R} A' \otimes_{\tilde A} \tilde \Gamma = R = A'[r]
\]
is defined by
\begin{align*}
g(a_2) &= a_2 + 3r\\
g(a_4) &= a_4 + 2 a_2 r + 3 r^2\\
g(a_6) &= a_4 r + a_2 r^2 + r^3
\end{align*}
This is a faithfully flat extension since $g(a_6)$ is a monic polynomial in $r$, and we invoke \fullref{basechg} to conclude that $H^*(A,\Gamma) \cong H^*(A',\Gamma')$ with
\[
(A',\Gamma') = \left(\Z_{(3)}[a_2,a_4], A'[r]/\left(r^3+a_2r^2+a_4r\right)\right).
\]
Let $I_0 = (3)$, $I_1 = (3,a_2)$, and $I_2 = (3,a_2,a_4)$ be the invariant prime ideals of $(A',\Gamma')$, and denote by $(A_n, \Gamma_n)$ the Hopf algebroid $(A',\Gamma') / I_n$ for $n=0,\;1,\;2$.

As an $A'$--module, $\Gamma'$ is free with basis $\{1,r,r^2\}$. We first compute the cohomology of $(A_2,\Gamma_2) = (\ZZ/3,\ZZ/3[r]/(r^3))$. Since $A_2$ is a field, $\Gamma_2$ is injective as an $(A_2,\Gamma_2)$--comodule, and we can take the minimal resolution
\[
\ZZ/3 \xrightarrow{\eta_L}  \Gamma_2 \xrightarrow{\frac\partial{\partial r}}  \Gamma_2 \xrightarrow{\frac{\partial^2}{\partial r^2}} \Gamma_2 \xrightarrow{\frac\partial{\partial r}} \Gamma_2 \xrightarrow{\frac{\partial^2}{\partial r^2}} \cdots
\]
and get the cohomology algebra
\[
H^{**}(A_2,\Gamma_2) = E(\alpha) \otimes P(\beta),
\]
where $\alpha \in H^{1,4}(A_2,\Gamma_2)$ and $\beta \in H^{2,12}(A_2,\Gamma_2)$.

In the cobar complex, $\alpha$ is represented by $[r]$ and $\beta$ by $[r^2|r]-[r|r^2]$. Furthermore, $\beta$ can be expressed as a Massey product:
\begin{equation}\label{masseybeta}
\beta = \langle \alpha, \alpha, \alpha \rangle.
\end{equation}
If we consider the algebraic Bockstein spectral sequence associated with the ideal $(a_4) \triangleleft A_1 = \Z/3[a_4]= A'/(3,a_2)$, we see that the polynomial generator corresponding to $a_4$ supports no differentials and hence
\[
H^{**}(A_1,\Gamma_1) = H^{**}(A_2,\Gamma_2) \otimes P(a_4).
\]
We now look at the Bockstein spectral sequence for $(a_2) \triangleleft A_0=\Z/3[a_2,a_4] = A'/(3)$. We find that in the cobar complex,
\begin{align*}
a_4 a_2^n &\mapsto 2a_2^{n+1}[r]\\
a_4[r] + a_2[r^2] & \mapsto 0 & \text {Massey product: } \langle a_2^{n+1}, \alpha,\alpha \rangle \\
a_4^2 a_2^n & \mapsto a_4 a_2^{n+1}[r]\\
a_4^2 [r] + a_4a_2[r^2] & \mapsto a_2^2 \left([r^2 | r] + [r|r^2]\right) = a_2^2 \beta\\
a_4^3 & \mapsto 0.
\end{align*}
Hence we get the following picture for the $E_\infty$ term:
$$\includegraphics{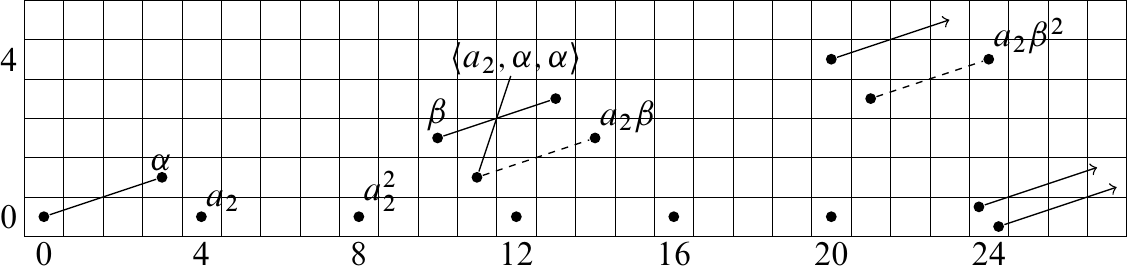}$$
Here, and in all subsequent charts, we use Adams indexing, ie, the square at coordinate $(x,y)$ represents the group $H^{y,x+y}$. We also use the following conventions:
\begin{itemize}
\item a dot `$\bullet$' represents a generator of $\Z/p$ (here, $p=3$);
\item a circle around an element denotes a nontrivial $\Z/p$--extension of the group represented by that element;
\item a square `$\square$' denotes a generator of $\Z_{(p)}$;
\item a line of slope $\infty$, $1$ or $1/3$ denotes multiplication with $p$, the generator of $H^{1,2}$ (not present at $p=3$), or of $H^{1,4}$;
\item an arrow with positive slope denotes generically that the pattern is continued in that direction;
\item an arrow with negative slope denotes a differential; the arrow, the source, and the target are then drawn in a gray color to denote that these classes do not survive.
\end{itemize}

The dashed line of slope $1/3$ denotes an exotic multiplicative extensions which follows from the Massey product shuffling lemma \cite[Appendix 1]{Ravenel:green2} and \eqref{masseybeta}:
\[
\langle a_2,\alpha,\alpha \rangle \alpha = a_2 \langle \alpha,\alpha,\alpha \rangle = a_2 \beta.
\]
Finally, we run the Bockstein spectral sequence to get the $3$--local homology. The 
differentials derive from $d_1(a_2) = 3r$ and yield
$$\includegraphics{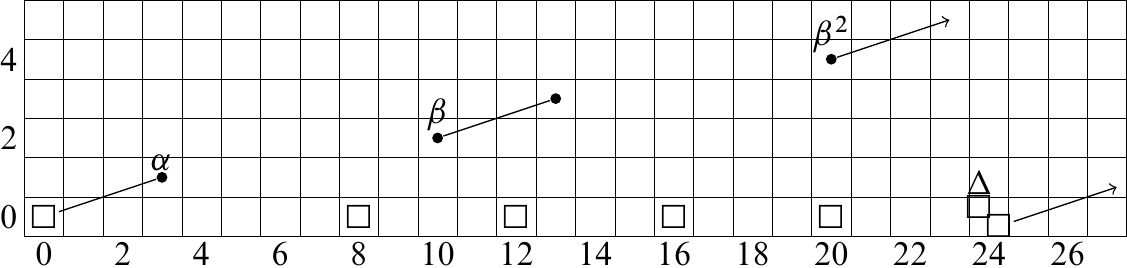}$$
where $\Delta$ is the discriminant.

\section{The differentials} \label{differentialsatthree}

In this section, I compute the differentials in the elliptic spectral sequence whose $E_2$--term was computed in the previous section. The differential which implies all further differentials in the $3$--local $\tmf$ 
spectral sequence is
\begin{equation}\label{dfiveatthree}
d_5(\Delta) = \beta^2 \alpha.
\end{equation}
This follows from comparison with the Adams--Novikov spectral sequence: By \eqref{ANSScomparison}, the class $[x_2]$ representing the Hopf map $\alpha_1 \in \pi_3(\SS^0)$ in the cobar complex of the Adams--Novikov spectral sequence is mapped to $\alpha$; because of the Massey product $\beta = \langle \alpha,\alpha,\alpha \rangle$, we also conclude that $\beta$ represents the Hurewicz image of the class $\beta_1 \in \pi_{10}(\S^0)$.

Now in the stable stems, $\beta_1^3 \alpha_1$ represents the trivial homotopy class although $\beta^3\alpha$ is nontrivial in the Hopf algebroid cohomology; thus it must be killed by a differential. The only possibility is $\beta^3\alpha = d_5(\Delta\beta)$, which implies the differential \eqref{dfiveatthree}.

\begin{figure}[hb]\small
\cl{\includegraphics{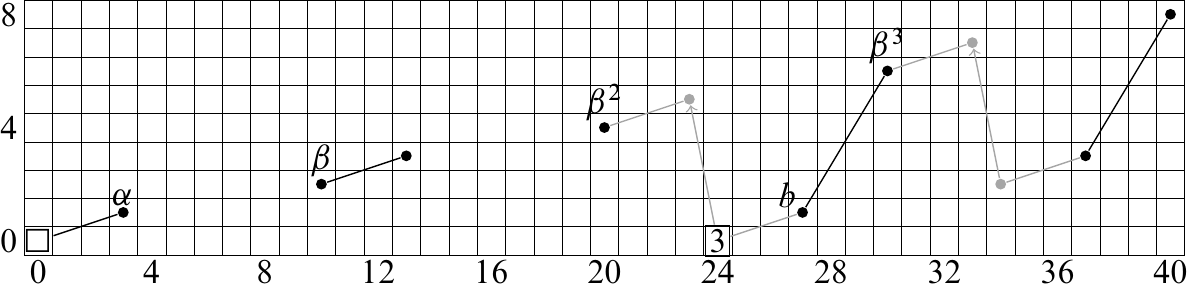}}
\caption{The elliptic spectral sequence at $p=3$ up to dimension $40$} \label{finalchart31}
\end{figure}

It also shows that there is an element $b$ representing a class in $\pi_{27}\tmf$ which is the unique element in the Toda bracket $\langle \beta_1^2,\alpha_1,\alpha_1\rangle$. The Massey product $\langle \alpha,\alpha,\alpha \rangle$ converges to the Toda bracket $\langle \alpha_1, \alpha_1, \alpha_1\rangle$. The Shuffling Lemma for Toda brackets \cite[Appendix 1]{Ravenel:green2} implies that there is a multiplicative extension $b \alpha_1 = \beta_1^3$ in $\pi_{27}$ as displayed in the chart. It implies that there has to be a differential $d_9(\Delta^2\alpha)=2\beta^5$ since $\beta_1^5 = \beta_1^2 b \alpha_1$.

\begin{figure}[ht]
$$\includegraphics{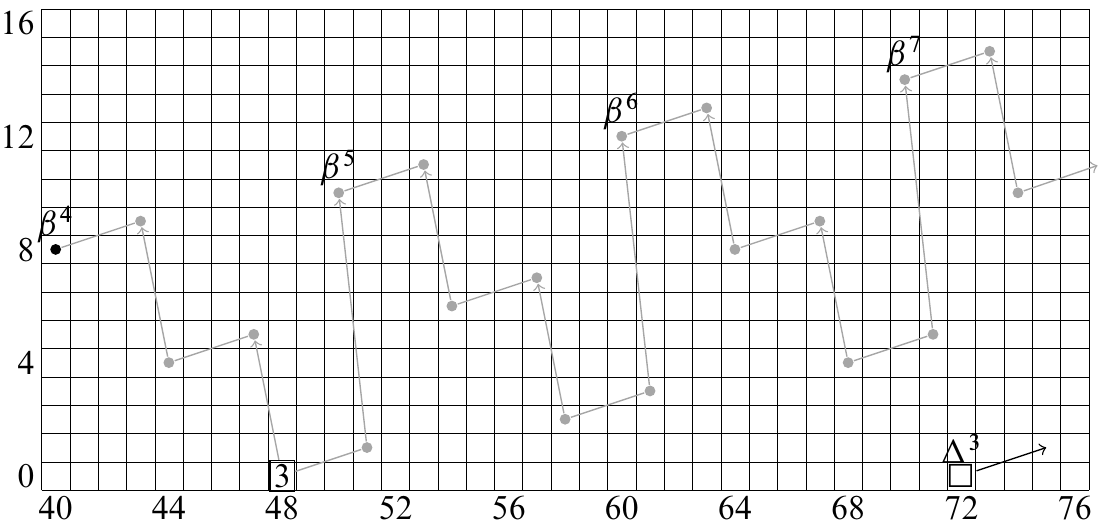}$$
\caption{The elliptic spectral sequence at $p=3$, dimensions 40--76} \label{finalchart32}
\end{figure}

The complete $E_\infty$ term has a free polynomial generator $\Delta^3$ in degree $72$ and is displayed in Figures \ref{finalchart31} and \ref{finalchart32}, where most classes in filtration $0$ have been omitted.

\section{The cohomology of the elliptic curve Hopf algebroid at\break \texorpdfstring{$p=2$}{p=2}} \label{exttwo}

In this section we are going to compute the cohomology of $(A,\Gamma)$, localized at the prime $2$.

Note that if $3$ is invertible, we can proceed as in \fullref{completethecube} and complete the cube in $x$ by a transformation $x \mapsto x+r$, thereby eliminating the coefficient $a_2$. This means that the resulting Hopf algebroid is
\[
(\tilde A,\tilde\Gamma) = (\ZZ_{(2)}[a_1,a_3,a_4,a_6], \tilde A[s,t]),
\]
where the map $\Gamma \to \tilde \Gamma$ maps $r$ to $\frac13(s^2+a_1 s)$.

By two applications of \fullref{basechg}, we can furthermore eliminate $a_4$ and $a_6$ ($a_4$ maps to a monic polynomial $s^4+$lower terms in $s$, $a_6$ to a monic polynomial $t^2$+lower terms in $t$). This introduces new relations in $\Gamma$:
\begin{gather*}
s^4 - 6 s t + a_1 s^3 - 3 a_1 t - 3 a_3 s = 0 \quad \text{and}\\
s^6 - 27 t^2 + 3a_1s^5 - 9a_1s^2t + 3a_1^2s^4 - 
  9a_1^2st + a_1^3s^3 - 27a_3t = 0
\end{gather*}

Hence, at the prime $2$, the Hopf algebroid 
\[
(A',\Gamma') = \left(\ZZ_{(2)}[a_1,a_3], A'[s,t]/(\text{the above relations})\right)
\]
has the same cohomology as $(A,\Gamma)$.

Note that as an $A'$--module, $\Gamma'$ is free with basis 
$\{1,s,s^2,s^3,t,st,s^2t,s^3t\}$ and we have invariant prime ideals $I_0 = (2)$, $I_1 = (2,a_1)$, and $I_2 = (2,a_1,a_3)$. Define $(A_n,\Gamma_n) = (A',\Gamma')/I_n$ for $n=0,\;1,\;2$.

We are going to compute the cohomology of this Hopf algebroid with a sequence of algebraic Bockstein spectral sequences. As in the $3$--primary case, we first compute $H^{**}(A_2,\Gamma_2)$. Observe that this is isomorphic as a Hopf algebra to the 
double $DA(1)$ of $A(1)$, where $A(1)$ is the dual of the subalgebra of the Steenrod algebra generated by $\Sq^1$ and $\Sq^2$. 

The cohomology of $DA(1)$ is well-known. I will compute it here nevertheless to determine cobar representatives of some classes.

Note that in $DA(1)$, $s^2$ is primitive, hence we can divide out by it to get a Hopf algebra which is isomorphic to the Hopf algebra $E(s) \otimes E(t)$. The Cartan--Eilenberg spectral sequence associated to the ideal $(s^2)$ coincides with the algebraic Bockstein spectral sequence of \fullref{algbocksteinss} in this situation and has
\[
E_1 = H^{**}(DA(1)/(s^2)) \otimes P(s^2) \Longrightarrow H^{**}(DA(1)).
\]
The cobar complex gives us the following differentials:
\begin{align*}
d_1[s^{2n}t]& = [s^{2n+2} | s];\\
d_1\left([t|t] + [s|s^2t] + [st|s^2]\right) &= [s^2 | s^2 | s^2].
\end{align*}
No more differentials are possible for dimension reasons, yielding the following $E_\infty$--term:
$$\includegraphics{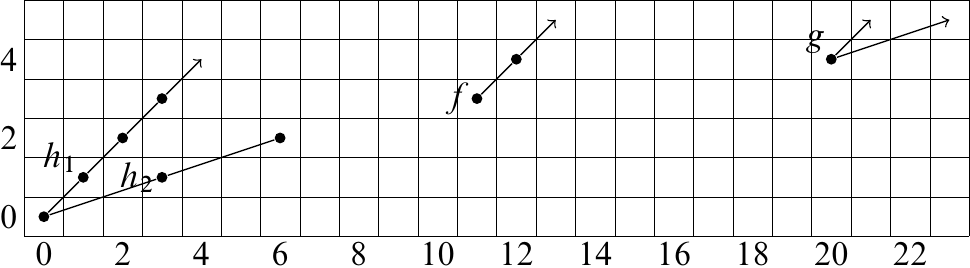}$$
The elements labeled $h_i$ are represented by $s^i$; the element $f$ is represented by
\[
f = [s|t|t] + [s|s|s^2t] + [s|st|s^2+t|s^2|s^2] = \langle h_2^2,h_2,h_1 \rangle,
\]
and the class $g$ is a polynomial generator in this algebra. It is represented in the cobar complex by the class $[t|t|t|t]$ plus elements of higher filtration in the above spectral sequence.
%
It can be expressed as a Massey product
\begin{equation} \label{firstmasseyg}
g = \langle h_2^2,h_2,h_2^2,h_2 \rangle.
\end{equation}
as can easily be seen from the differentials.

We record the following Massey product relation:
\begin{equation} \label{masseyh22}
h_2^2 = \langle h_1,h_2,h_1 \rangle
\end{equation}
which follows from $d[t] = [s^2 | s]$, $d\left([t]+[s^3]\right) = [s | s^2]$, and
\[
d[ts] = [t|s] + [s|t+s^3] + [s^2 | s^2].
\]
We will now start an algebraic Bockstein spectral sequence to compute the cohomology of $(A_1,\Gamma_1)$ from $(A_2,\Gamma_2)$:
\[
E_1 = H^{**}(DA(1)) \otimes P(a_3) \Longrightarrow H^{**}(A_1,\Gamma_1).
\]
$$\includegraphics{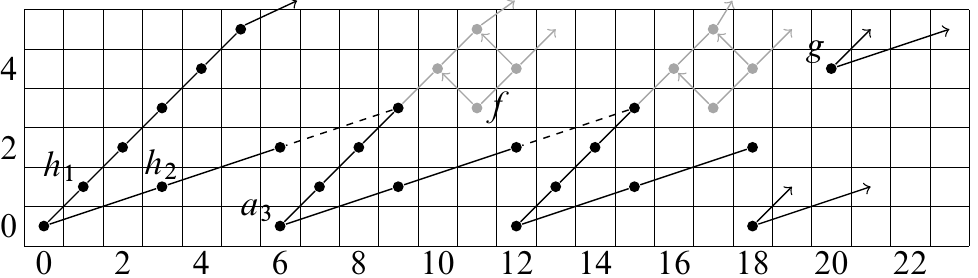}$$
The elements $h_1$ and $h_2$ lift to cycles in the $E_\infty$ term, but the cobar complex shows that
\[
d_1(f) = a_3 h_1^4.
\]
We have an extension $a_3 h_1^3 = h_2^3$ because in the cobar complex,
\[
[t | t] + [s | s^2t] + [st | s^2] \mapsto [s^2 | s^2 | s^2] + a_3 [s|s|s],
\]
as well as an exotic Massey product extension
\begin{equation}\label{masseyc}
a_3h_1^2 = \langle h_2,h_1,h_2 \rangle.
\end{equation}
Indeed, the Massey product has the cobar representative
\[
[t | s^2] + [s^2 | t+s^3],
\]
and
\[
d[s^2(t+s^3)] = [t|s^2] + [s^2|t+s^3] + [s^4|s], \text{ and } s^4 = a_3 s \pmod{2,a_1}.
\]
No more differentials are possible for dimension reasons.

We thus have computed
\begin{align*}
E_\infty &= P(a_3,h_1,h_2,g)/(h_1h_2,h_2^3+a_3h_1^3)\\
\intertext{or, additively only,}
E_\infty&= P(a_3,g)\{h_2,h_2^2,h_1^i \mid i \geq 0\}/(a_3h_1^3)
\end{align*}
It will be necessary to keep track of a Massey product expression for $g$. Equation~\eqref{firstmasseyg} becomes invalidated because $h_2^3 \neq 0$; instead, \eqref{firstmasseyg} lifts to a matric Massey product representation
\begin{equation}\label{secondmasseyg}
g = \fourmatricmassey{h_2^2}{h_1}{h_2}{a_3 h_1^2}.
\end{equation}
The next step is to introduce $a_1$. The $d_1$--differential is linear over $P(a_3^2,g)$ and satisfies
\[
d_1(a_3 h_2^i) = a_1h_2^{i+1}
\]
where $h_2$ is the lift of the cycle previously called $h_2$, represented by $[s^2+a_1s]$ in the cobar complex.

The resulting homology may be written additively as
\[
P(g) \otimes \left( K \oplus P(a_1) \otimes L \oplus H_1\right)
\]
where
\begin{align*}
K &= P(a_3^2)\left\{h_2,h_2^2,h_2^3\right\}
\intertext{is the $a_1$--torsion,}
L &= P(a_3^2)\left\{h_1^i, x, xh_1 \mid 0 \leq i \leq 3\right\},
\intertext{and}
H_1 &= \F_2\left\{h_1^i \mid i \geq 4\right\}
\end{align*}
is the initial infinite $h_1$--tower. Here
\begin{equation} \label{masseyx}
x = \langle a_1,h_2,h_1 \rangle
\end{equation}
is cobar represented by $[a_3s+a_1t]$.

The submodule $K \oplus 1 \otimes L+H_1$ constitutes the filtration $0$ part; this is displayed in the following diagram, along with the $h_1$-- and $h_2$--multiplications:
$$\includegraphics{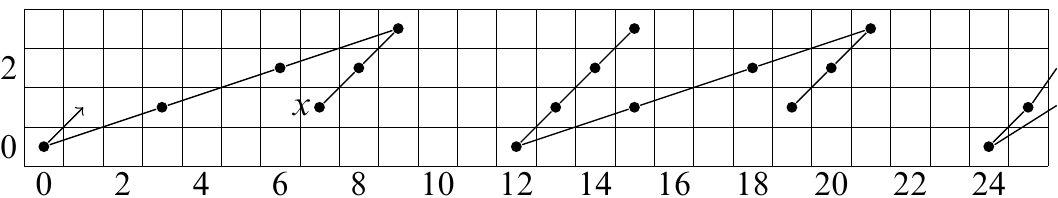}$$
The $d_2$--differentials all follow from $d_2(a_3^2)=a_1^2x$; explicitly,
\[
d_2(a_3^{2k} a_1^n) = \begin{cases} a_3^{2k-2} a_1^{n+2} x; & k \text{ odd};\\
0; & k \text{ even.}\end{cases}
\]
Since $d_2$ does not involve $a_1$--torsion and is $(g,a_3^4)$--linear, it is enough to understand what it does on $L/(a_3^4,g)$. This module decomposes into two $d_2$--invariant copies $L' \oplus a_1L'$ of a module $L'$, which is shown in the following chart:
$$\includegraphics{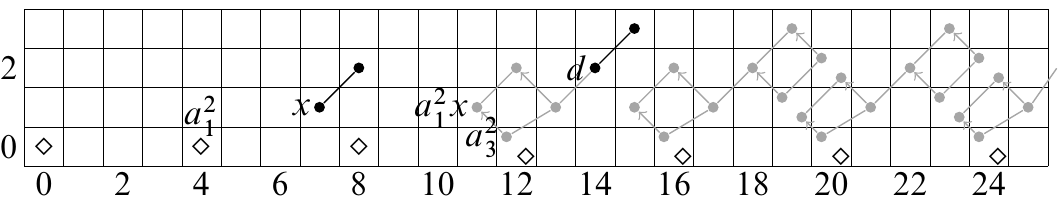}$$
Diamonds denote $h_1$--strings of $P(h_1)/(h_1^4)$ of length 4. Thus
\[
H_*(L',d_2) =  P(a_1^2,h_1)/(h_1^4) \oplus \F_2\left\{x,xh_1,d,dh_1\right\} \quad \text{where } d=x^2.
\]
The next chart displays the full spectral sequence modulo Bockstein filtration~5 (the full chart would become too cluttered). Gray classes that do not seem to be hit by or to support a differential support a differential into higher filtration.
$$\includegraphics{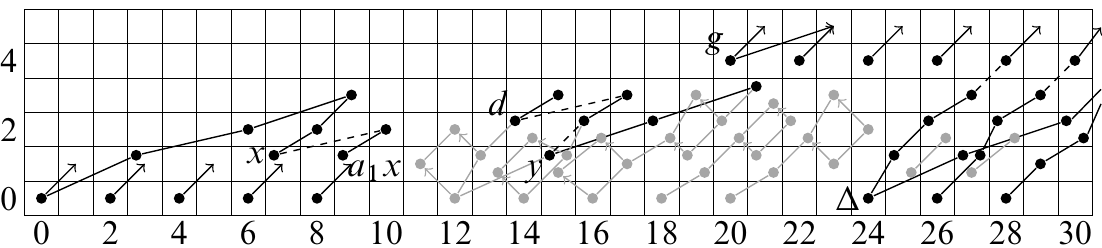}$$
The new indecomposable class $y$ has the Massey product representation:
\begin{equation}\label{masseyy}
y = \langle x, a_1^2,h_2 \rangle = \langle a_1x,a_1,h_2 \rangle,
\end{equation}
and the class $\Delta$ can be expressed as
\begin{equation}\label{masseydelta}
\Delta \in \langle a_1x,a_1,a_1x,a_1 \rangle \quad \text{with indeterminacy } \subseteq a_1 H^{0,22}(A_0,\Gamma_0)
\end{equation}
There are a number of exotic multiplicative extensions, which I will now verify.

\begin{lemma} \label{a1extensions}
The following multiplicative extensions exist in the above algebraic Bockstein spectral sequence:
\begin{enumerate}
\item $x h_2 = a_1 x h_1$ \label{a1ext1}
\item $d h_2 = a_1 d h_1$ \label{a1ext2}
\item $y h_1 = a_1 d$ \label{a1ext3}
\item $a_1^4 g = \Delta h_1^4$ \label{ga14extension}
\end{enumerate}
\end{lemma}
\begin{proof}
\eqref{a1ext1}\qua
This follows from the Massey product Shuffling Lemma (cf \cite[Appendix 1]{Ravenel:green2}), $\langle a,b,c \rangle d = \pm a \langle b,c,d \rangle$:
\[
x h_2 \underset{\eqref{masseyx}}= \langle a_1,h_2,h_1 \rangle h_2 = a_1 \langle h_2,h_1,h_2 \rangle \underset{\eqref{masseyc}}= a_1 x h_1.
\]

\eqref{a1ext2}\qua This follows from \eqref{a1ext1} by multiplication with $x$ since $d=x^2$.

\eqref{a1ext3}\qua \(\displaystyle{
y h_1 \underset{\eqref{masseyy}}= \langle x,a_1^2,h_2\rangle h_1 = x \langle a_1^2,h_2,h_1 \rangle = a_1 x^2 = a_1 d.
}\)

\eqref{ga14extension}\qua This is the hardest extension to verify. We first claim:
\begin{equation}\label{ga12expression}
g a_1^2 = \langle h_1,dh_1,x \rangle.
\end{equation}
This implies \eqref{ga14extension} since then,
\[
g a_1^4 = \langle h_1,dh_1,x \rangle a_1^2 = h_1 \langle dh_1,x,a_1^2 \rangle.
\]
Thus $g a_1^4$ is divisible by $h_1$, and the only possibility is $g a_1^4 = \Delta h_1^4$.

To prove \eqref{ga12expression}, we start from the Massey product expression for $g$ and compute
\begin{multline*}
ga_1^2 \underset{\eqref{secondmasseyg}}= \fourmatricmassey{h_2^2}{h_1}{h_2}{xh_1}a_1^2
\\
\subseteq \left\langle \begin{pmatrix}h_2^2 & h_1\end{pmatrix},
\begin{pmatrix}0 & a_1 x h_1 \\ a_1 x h_1 & 0\end{pmatrix},
\begin{pmatrix}0 & a_1 h_1 \\ a_1 h_1 & 0\end{pmatrix},
\begin{pmatrix}h_2\\x h_1\end{pmatrix} \right \rangle\\
= \langle h_1,a_1xh_1,a_1h_1,xh_1 \rangle + \langle h_2^2,a_1xh_1,a_1h_1,h_2 \rangle = L+R.
\end{multline*}
The second summand, $R$, is not strictly defined in the technical sense that the sub-Massey products have indeterminacy, as
\[
\langle a_1xh_1,a_1h_1,h_2\rangle = \{ 0, dh_2 \}.
\]
Despite this, the four-fold product $R$ has no indeterminacy and contains only $0$. To see this, we calculate
\[
a_1 R \subseteq \langle a_1,h_2^2,a_1xh_1,a_1h_1 \rangle h_2,
\]
where the Massey product is now strictly defined, showing that $a_1 R$ is divisible by $h_2$, which implies that it is $0$. Multiplication with $a_1$ in that degree is injective, hence $R=0$.

As for the first summand, we find that
\[
L = \langle h_1,xh_2,h_1,xh_2\rangle\\
\underset{\text{\fullref{shufflinglemma}}}= \langle h_1,\langle xh_2,h_1,h_2\rangle,x\rangle = \langle h_1,dh_1,x \rangle.
\]
Equality holds because both left and right hand side have zero indeterminacy.
\end{proof}

\begin{remark}
One might be tempted to think that $\Delta h_1^4 = x^4$, but this cannot hold here since $\Delta h_1^8 \neq 0$ whereas $(xh_1)^4=0$ because $xh_1^3 = 0$.
\end{remark}
 
For proving \ref{a1extensions}\eqref{ga14extension}, we made use of the following lemma:
\begin{lemma} \label{shufflinglemma}
Let $a \in M$, $b$, $c$, $d$, $e \in A$ for a differential graded module $M$ over a differential graded $\F_2$--algebra $A$. Suppose $ab=bc=cd=0$ and $\langle a,b,c \rangle = 0$. Then
\begin{equation*}
\langle a, \langle b,c,d \rangle, e \rangle \cap \langle a,b,c,d\,e \rangle \neq \emptyset.
\end{equation*}
\end{lemma}

\begin{proof}
We adopt the following notation: for a boundary $x$, we denote by $\undl{x}$ a chosen chain such that $d(\undl x)=x$, keeping in mind that it is not unique. Consider the following defining system for $\langle a,b,c,d\,e \rangle$:
\begin{equation*}
\xymatrix@!0{
a && b && c && d\,e\\
& {\undl{ab}} && {\undl{bc}} && {\undl{cd}\,e}\\
&& {\undl{\langle a,b,c \rangle}} && {\undl{\langle b,c,d\,e \rangle}}
}
\end{equation*}
Note that this is not the most general defining system because we insist that the class bounding $c\,d\,e$ actually is a class $\undl{cd}$ bounding $c\,d$, multiplied with $e$.

On the other hand, a defining system for $\langle a, \langle b,c,d \rangle, e \rangle$ is given by
\begin{equation*}
\xymatrix@!0{
a && {\undl{bc}d+b\undl{cd}} && e\\
& {\begin{matrix}
\undl{\langle a,b,c \rangle}d\\
+ \undl{ab}\,\undl{cd}
\end{matrix}
}
&& {\undl{\langle b,c,d\,e\rangle}}.
}
\end{equation*}
If we compute the representatives of the Massey product for both defining systems, we get in both cases:
$$a \undl{\langle b,c,d\,e \rangle} + \undl{\langle a,b,c \rangle} d\,e + \undl{a\,b}\, \undl{c\,d}\, e.\proved$$
\end{proof}

The last step is to run the Bockstein spectral sequence to compute the integral cohomology. We have that 
\begin{align*}
d_1(a_1^{2k+1}) &= 2a_1^{2k} h_1\\
d_2(a_1^2) &= 4 h_2\\
d_4(yh_2^2) &= 8g
\end{align*}
The first two are immediate, and the last one follows once again from a Massey product expression for $g$ in the $E_4$ term of the BSS.

Setting $c=\langle h_2,h_1,h_2 \rangle$ (cf \eqref{masseyc}), we have
\[
g = \fourmatricmassey{h_2^2}{h_1}{h_2}{c}.
\]
Note that
\[
d_1(x) = d_1\langle a_1,h_2,h_1 \rangle = 2 \langle h_1,h_2,h_1 \rangle \underset{\eqref{masseyh22}}= 2h_2^2
\]
Thus, when multiplying $g$ by $2$ from the left,
\begin{align*}
\left\langle 2,\begin{pmatrix}h_2^2 & h_1\end{pmatrix},
\begin{pmatrix}h_2\\ c\end{pmatrix}\right\rangle &= 0 \quad \text{and}\\
\left\langle 2,\begin{pmatrix}h_2^2 & h_1\end{pmatrix},
\begin{pmatrix}c\\h_2\end{pmatrix}\right\rangle &= \langle 2,h_2^2,c \rangle = dh_1.
\end{align*}
Hence
\begin{equation}\label{massey2g}
2g = \left\langle dh_1,\begin{pmatrix} h_1 & h_2^2 \end{pmatrix},\begin{pmatrix} h_2\\c \end{pmatrix} \right\rangle = \langle dh_1,h_1,h_2 \rangle.
\end{equation}
From this it follows that
\begin{equation}\label{gextension}
4g = 2\langle dh_1,h_1,h_2 \rangle = \langle 2,dh_1,h_1 \rangle h_2 = dh_2^2.
\end{equation}
Therefore, the differential follows from this exotic $h_2$ extension and $d_1(y) = 2d$. This in turn follows by multiplication with $h_1$: $d_1(yh_1) = d_1(a_1d) = 2dh_1$.

Let us also record the identity
\begin{equation}\label{masseyd}
d = \langle h_2,2h_2,h_2,2h_2 \rangle.
\end{equation}
The cohomology of $(A',\Gamma')$ is given by the page below:
$$\includegraphics{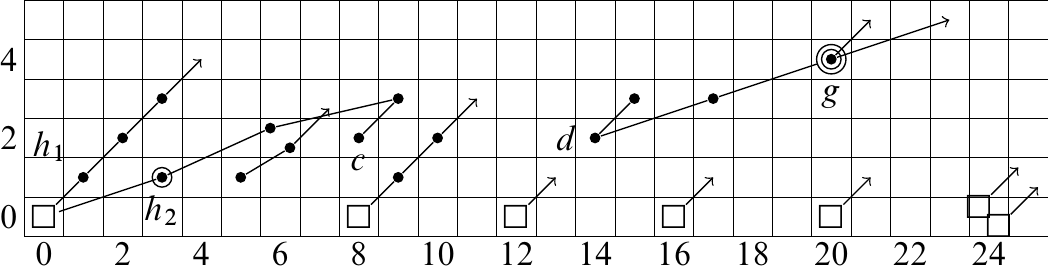}$$

\section{The differentials} \label{differentialsattwo}

In this section, I compute the differentials in the spectral sequence whose $E_2$--term was computed in the previous section. The first couple of differentials can be derived by comparing this spectral sequence to the Adams--Novikov spectral sequence; all the further differentials follow from Toda bracket relations; in fact, all those differentials are determined from knowing that $\tmf$ is an $A_\infty$--ring spectrum.

Under the Hurewicz map $H$ of \fullref{ellalg}, the classes I called
$h_1$ and $h_2$ in the previous section are the images of the classes
of the same name in the ANSS by \eqref{ANSScomparison}; in particular,
they will represent the Hurewicz images of $\eta \in \pi_1$ and $\nu
\in \pi_3$. Furthermore, as the map $\SS^0 \to \tmf$ is an
$A_\infty$--map, Toda brackets are mapped to Toda brackets, and Massey
products in the ANSS are mapped to Massey products in the elliptic
spectral sequence. By the computations of Massey products summarized
in \fullref{todaapp}, it follows that the classes $c$ and $d$ are the
images of the classes called $c_0$ and $d_0$ in the ANSS and represent
the homotopy classes $\epsilon \in \pi_8$ and $\kappa \in \pi_{14}$,
respectively.

Furthermore, the generator $g \in H^{4,24}(A',\Gamma')$ represents
\begin{equation}\label{kbartoda}
\copy \kbar = \langle \kappa,2,\eta,\nu \rangle \in \pi_{20} \tmf,
\end{equation}
the image of the class of the same name in $\pi_{20} \SS^0$. This
point requires a bit of explanation since the Toda bracket is not
reflected by a Massey product in the elliptic spectral sequence, and
because the corresponding class in the ANSS has filtration $2$ instead
of $4$. The given Toda bracket is computed by Kochman
\cite[5.3.8]{Kochman:computer} for the sphere.  Since both in the
sphere (by an exotic extension in the Adams--Novikov $E_\infty$--term)
and in $\tmf$ (by \eqref{gextension}), we have that $4\copy \kbar =
\kappa \nu^2$, we have that $\copy \kbar \in \pi_{20} \tmf$ must be an
odd multiple (mod 8) of $\copy \kbar \in \pi_{20} \SS^0$. We can alter
$g$ by an odd multiple to make sure they agree on the nose.

\subsection {\texorpdfstring{$d_3$}{d3}--differentials}
It is easy to see that there is no room in the spectral sequence for $d_2$--differentials (indeed, the checkerboard phenomenon implies that there are no $d_i$ for any even $i$). We want to describe all $d_3$ differentials.

Since $\eta^4$ is zero in $\pi_4^s$, $h_1^4$ has to be 
hit by a differential at some point. The only possibility for this is
\[
d_3(a_1^2 h_1) = h_1^4.
\]
This implies
\[
d_3(a_1^{4n+2}) = a_1^{4n} h_1^3.
\]
There is no more room in the spectral sequence for any differentials up to dimension $15$. We may thus study the multiplicative extensions in $E_5 = E_\infty$ in this range. The Hurewicz map dictates that there must be a multiplicative extension $4 h_2 = h_1^3$.

By multiplicativity, there are no more possibilities for a nontrivial $d_3$.

The following chart shows the $d_3$ differentials; black classes survive.
$$\includegraphics{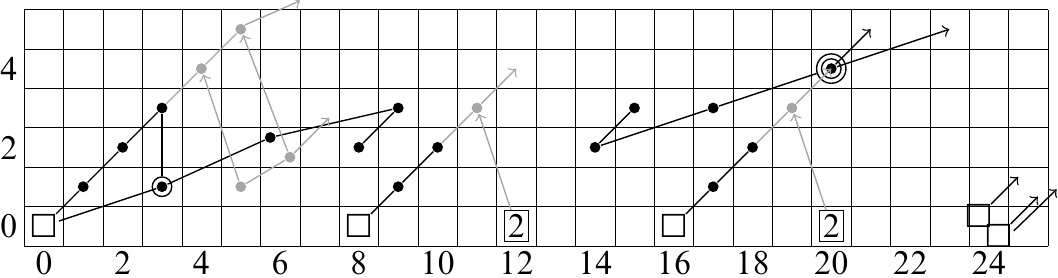}$$

\subsection{Higher differentials and multiplicative extensions}
For the determination of the higher differentials, we will only need to appeal once more to an argument involving the Hurewicz map from the sphere; the remaining differentials and exotic extension all follow algebraically using Massey products and Toda brackets. The crucial differential is
\[
d_5(\Delta) = g h_2.
\]
There are several ways to see this. Maybe the easiest is to notice that
\[
\copy \kbar \nu^3 = 0 \in \pi_{29}(\S^0)
\]
thus the class $g h_2^3$ representing its Hurewicz image has to be
killed by a differential. There is only one class in $E_5^{30+*,*}$,
namely, $\Delta h_2^2 \in E_5^{32,2}$. But if $d_5(\Delta h_2^2) = g
h_2^3$, then necessarily $d_5(\Delta) = g h_2$. Because of the
multiplicative extension $4 \nu = \eta^3$, there also has to be a
$d_7(4\Delta) = g h_1^3$.

The following chart displays the resulting $d_5$ and $d_7$
differentials up to dimension 52 along with the multiplicative
extensions they incur. From here on I omit the classes in filtration
$\leq 2$ which look like the Adams--Novikov chart of $bo$; it is a
straightforward verification that they can never support a
differential from $E_5$ on.
$$\includegraphics{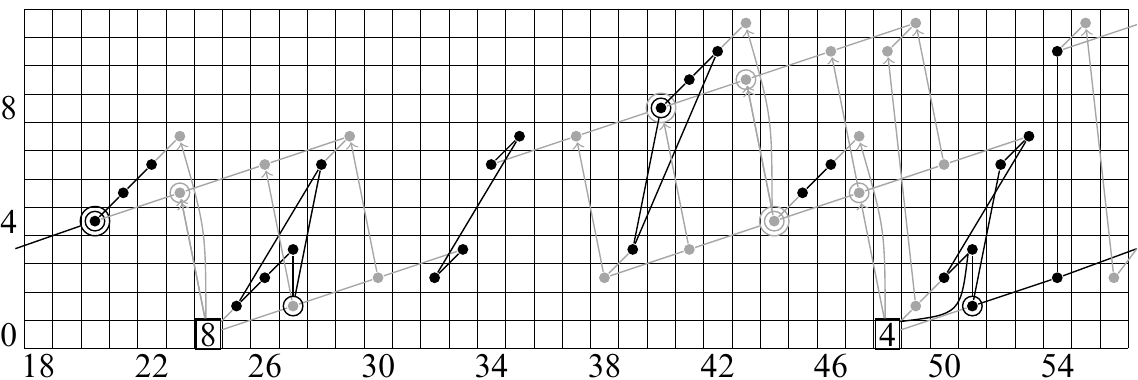}$$
Since $d_5(\Delta)=g h_2$, we have by multiplicativity that $d_5(\Delta^i) = i \Delta^{i-1} g h_2$. Thus the following powers are cycles for $d_5$: $\{\Delta^{4i},2\Delta^{4i+2},4\Delta^{2i+1}\}$. Since $d_7(4\Delta) = g h_1^3$, we furthermore have that the following classes survive to $E_9$:
\[
\{\Delta^{8i},2\Delta^{8i+4},4\Delta^{4i+2},8\Delta^{2i+1}\}.
\]
It will turn out that these are actually also the minimal multiples of the powers of $\Delta$ that survive to $E_\infty$. Thus, $\Delta^8$ is a polynomial generator in $\pi_{192}\tmf$.

I will now verify the extensions displayed in the above chart.

\begin{prop} \label{efiveextensions}
The following exotic multiplicative extensions exist in $\pi_*\tmf$. We denote a class in $\pi_{t}\tmf$ represented by a cycle $c \in E_5^{t+s,s}$ by $e[s,t]$ if $c$ is the only class in that bidegree that is not divisible by $2$.
\begin{enumerate}
\item $e[25,1] = \langle \copy \kbar,\nu,\eta\rangle$, thus $e[25,1]\nu = \copy \kbar \epsilon$ \label{firstnuextension}
\item $e[27,1] = \langle \copy \kbar,\nu,2\nu \rangle$, thus $e[27,1] \eta = \copy \kbar \epsilon$ \label{firstetaextension}
\item $e[32,2] = \langle \copy \kbar,\nu,\epsilon\rangle$, thus $e[32,2] \nu = \copy \kbar \kappa \eta$
\item $e[39,3] = \kappa e[25,1]$, thus $e[39,3] \eta = 2 \copy \kbar^2$ \label{extensionthirtynineeta}
\item $\copy \kbar \eta^2 = \kappa \epsilon$, thus $e[39,3] \nu = \copy \kbar^2 \eta^2$ \label{epsilonextension}
\end{enumerate}
\end{prop}

\begin{proof}
The representations as Toda brackets follow immediately from Massey products in the dga $(E_5,d_5)$.
For all multiplicative extensions, use the shuffling lemma for Toda brackets together with the identities from \fullref{todaapp}.

I will verify the interesting $\epsilon$ extension on $\kappa$. Since $\kappa \nu^2 = 4\copy \kbar$ and $4 \nu = \eta^3$, we find that $\kappa \epsilon \eta = \copy \kbar \eta^3$. Since the group in this degree is uniquely divisible by $\eta$, the result follows.
\end{proof}

This concludes the computation of the $E_9$ page. For degree reasons, we have that $E_9^{s,t} = E_\infty^{s,t}$ for $t-s \leq 48$; the first possible $d_9$ differential is
\[
d_9(e[49,1]) \overset{?}= g^2 c
\]
This differential is in fact present. In fact, if it were not, we would have an $\eta$ extension analogous to \ref{efiveextensions}\eqref{firstetaextension} from $2e[47,5]$ to $\copy \kbar^2 \epsilon$; but the representative of the former class is hit by the differential $d_5$ on $\Delta^2$. This propagates to the following picture of differentials:
$$\includegraphics{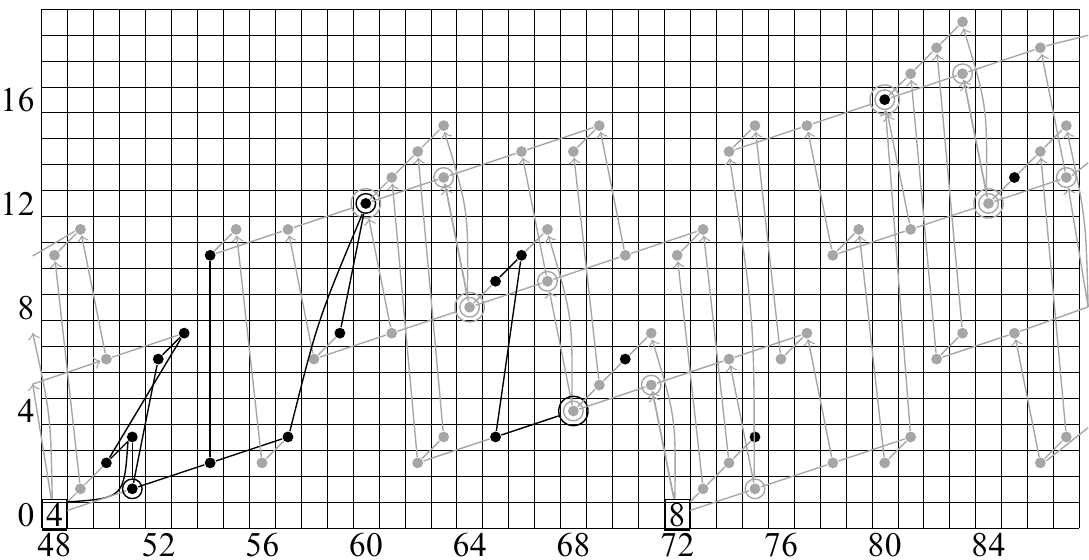}$$
I will first explain the induced $d_9$--differentials. The class $e[54,2]$ is in the Massey product $\langle h_2,e[49,1],h_2 \rangle$, and thus
\[
d_9(e[56,2]) = d_9(\langle  h_2,e[49,1],h_2 \rangle) = \langle h_2,g^2c,h_2\rangle = g^2dh_1.
\]
Secondly, the class $e[73,1]$ can be written as $\langle g,h_2,e[49,1]\rangle$; thus,
\[
d_9(e[73,1]) = \langle g,h_2,g^2 c \rangle = e[72,10].
\]
The remaining $d_9$--differentials follow easily from the multiplicative structure.

We now turn to the multiplicative extensions.

\begin{prop} \label{enineextensions} $\;$

\begin{enumerate}
\item $e[50,2] = \left\langle \copy \kbar^2,\begin{pmatrix} \epsilon & \nu \end{pmatrix},\begin{pmatrix}\eta\\ \nu^2 \end{pmatrix}\right\rangle$,
thus $e[50,2] \nu = e[53,7]$.
\item $e[51,1] \eta = e[52,6]$. \label{extensionfiftyone}
\item $e[54,2] = \langle \copy \kbar^2,\nu,2\nu,\nu\rangle\nu$, thus $2 e[54,2] = e[54,10]$. \label{extensionfiftyfour}
\item $e[57,3] = \langle \copy \kbar^2 ,\kappa \eta,\eta \rangle$, thus $e[57,3] \nu = 2\copy\kbar^3$.
\item $e[65,3] \eta = e[64,10]$.
\end{enumerate}
\end{prop}
\begin{proof}
\begin{enumerate}
\item The extension follows from the Toda bracket expression as follows:
\[
\left\langle \copy \kbar^2,\begin{pmatrix} \epsilon & \nu \end{pmatrix},\begin{pmatrix}\eta\\ \nu^2 \end{pmatrix}\right\rangle \nu = \langle \copy \kbar^2,\nu,\nu^3 \rangle = e[53,7].
\]
\item This is the same as \ref{efiveextensions}\eqref{firstetaextension}.
\item This follows from 
\[
\langle \copy\kbar^2,\nu,2\nu,\nu\rangle \cdot 2 \nu = \copy \kbar^2 \langle \nu,2\nu,\nu,2\nu \rangle = \copy  \kbar^2 \kappa.
\]
\item This follows from
\[
2\copy \kbar^3 = \copy\kbar^2 \langle \kappa \eta,\eta,\nu \rangle = \langle \copy \kbar^2 ,\kappa \eta,\eta \rangle \nu.
\]
\item This is the extension \eqref{extensionfiftyone} multiplied with $\kappa$, using \ref{efiveextensions}\eqref{epsilonextension}.\proved
\end{enumerate}
\end{proof}

The first higher differential occurs in dimension $62$, where
\[
d_{11}(e[62,2]) = g^3h_1.
\]
It suffices to show that $d_{11}(e[63,3]) = g^3h_1^2$. The latter class has to be hit by some differential since it is $gh_2$ times $e[39,3]$ by \fullref{efiveextensions}\eqref{epsilonextension}, and the claimed $d_{11}$ is the only possibility.

The $2$--extension in the $54$--stem also forces a
\[
d_{11}(e[75,2]) = g^3 d
\]
since the target class would have to be $2 e[74,6]$, which is a $d_5$--boundary. 

\begin{figure}[ht]
\bigskip
\cl{\includegraphics{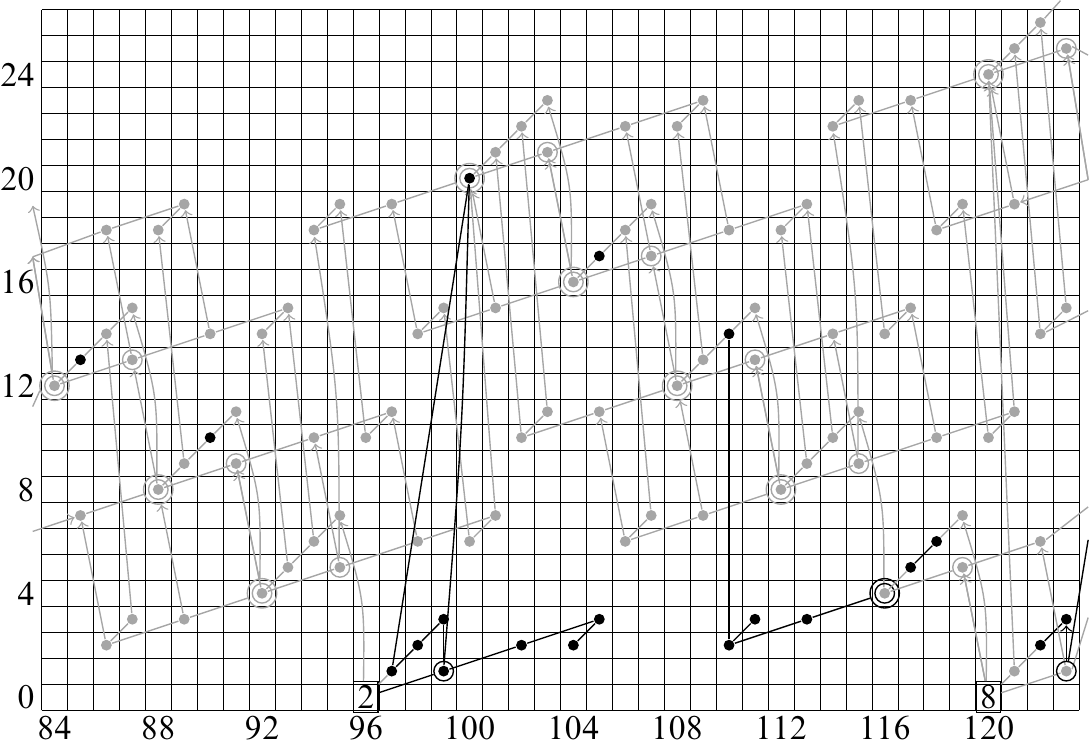}}
\caption{The elliptic spectral sequence --- dimensions 84--123} \label{chart84}
\end{figure}

In a similar fashion, the $\eta$--extension from \ref{efiveextensions}\eqref{extensionthirtynineeta} forces
\[
d_{13}(e[81,3]) = 2g^4.
\]
We proceed to the next chunk, starting in dimension $84$ and displayed in \fullref{chart84}.

Up to total dimension $119$, $E_{13} = E_\infty$ because there is no room for longer differentials, and the remaining $d_{11}$ and $d_{13}$ differentials follow from
the ones already computed by multiplicativity. There are some multiplicative extensions which I collect in the following proposition. The crucial part here is to express $\copy \kbar$ as a Toda bracket in a different fashion.

\begin{lemma} \label{kbarlemma}
We have
\[
\copy \kbar = \langle \nu,2\nu,\nu,4\nu,\nu,\eta\rangle = \langle \nu,2\nu,\nu,4\nu,\eta,\nu\rangle.
\]
\end{lemma}
\begin{proof}
We start with the representation $\copy \kbar = \langle \kappa,2,\eta,\nu\rangle$ (cf \fullref{todaapp}). Then
\begin{multline*}
\eta \copy \kbar = \langle \kappa \eta,2,\eta,\nu \rangle = \langle \langle \epsilon,\nu,2\nu\rangle,2,\eta,\nu \rangle \\
= \langle \epsilon,\nu,4\nu,\eta,\nu \rangle = \langle \langle \eta,\nu,2\nu \rangle,\nu,4\nu,\eta,\nu \rangle
= \eta \langle \nu,2\nu,\nu,4\nu,\eta,\nu \rangle.
\end{multline*}
Since both classes are uniquely divisible by $\eta$, the first part follows. For the second equality, the right hand side can be split up as
$$
\langle\langle \nu,2\nu,\nu,2\nu\rangle,2,\eta,\nu\rangle = \langle \kappa,2,\eta,\nu \rangle.\proved$$
\end{proof}

\begin{corollary} \label{ethirteenextensions}\ 

\begin{enumerate}
\item $e[97,1] = \langle \copy \kbar^4,\nu,2\nu,3\nu,4\nu,\eta\rangle$, thus $e[97,1]\nu = \copy \kbar^5$.
\item $e[99,1] = \langle \copy \kbar^4,\nu,2\nu,3\nu,4\nu,\nu\rangle$, thus $e[99,1]\eta = \copy \kbar^5$.\label{extensionninetynine}
\item $2 e[110,2] = e[110,14]$
\end{enumerate}
\end{corollary}
\begin{proof}
The first two parts are immediate using the above lemma.

The last part requires further explanation. We have
\[
e[110,2] = \langle \copy \kbar^4,\nu,2\nu,\nu,4\nu,\kappa \rangle
\]
and thus
$$
2e[110,2] = \langle \copy \kbar^4,\nu,2\nu,\langle \nu,4\nu,\kappa,2 \rangle\rangle = \langle \copy \kbar^4,\nu,2\nu,\langle 2,\copy\kbar,2 \rangle \rangle$$
$$
= \langle \copy \kbar^4,\nu,2\nu,\copy\kbar \eta \rangle 
= \langle \copy \kbar^5,\nu,2\nu,\eta^2\rangle = e[110,14].\proved$$
\end{proof}

There is one long differential
\[
d_{23}(e[121,1]) = g^6
\]
which follows immediately from the $\eta$--extension \ref{ethirteenextensions}\eqref{extensionninetynine}.

\begin{figure}[ht]
\cl{\includegraphics{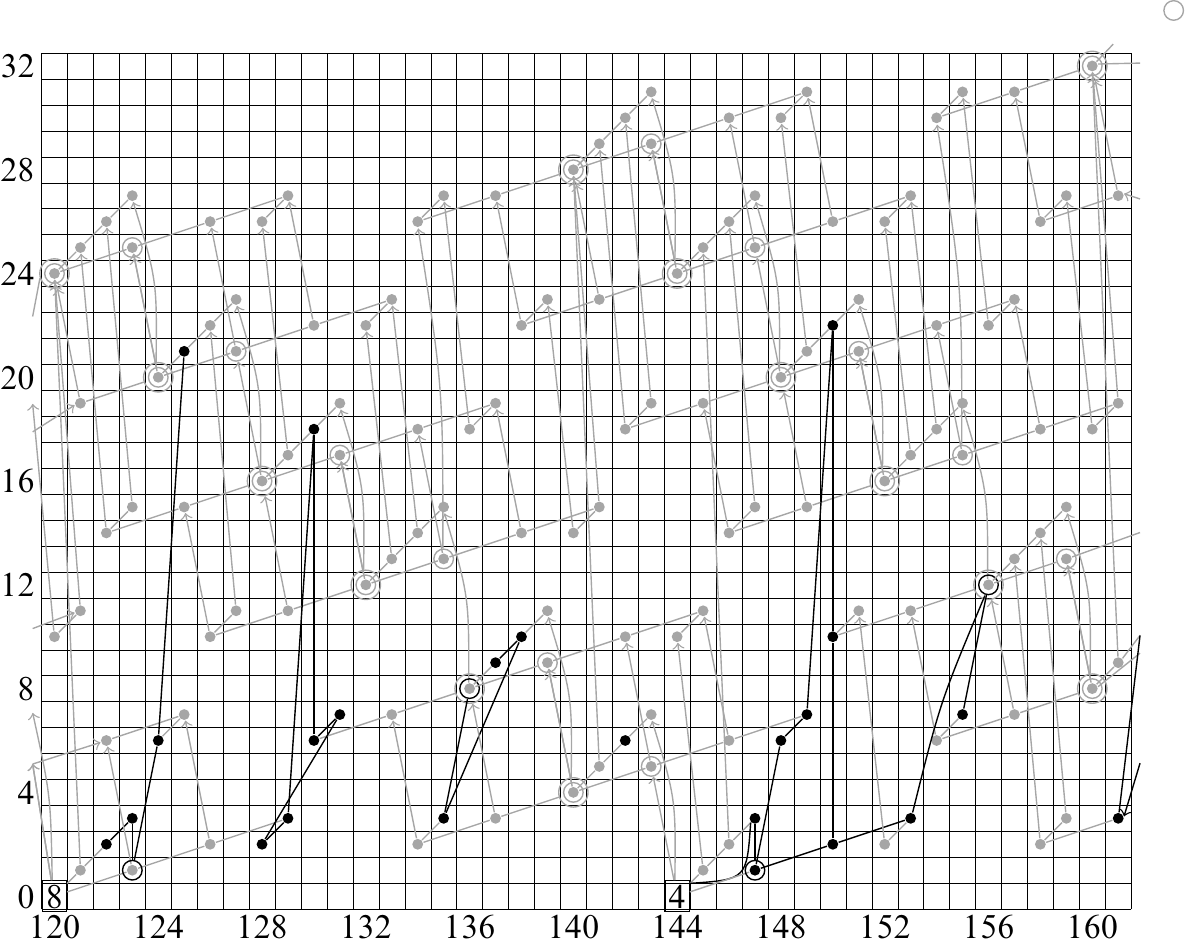}}
\caption{The elliptic spectral sequence --- dimensions 120--161}\label{chart120}
\end{figure}

Far from running out of steam, we press on to the next $40$ dimensions, as displayed in \fullref{chart120}.

Again, we address the multiplicative extensions indicated in the chart. The extensions that raise filtration by at most nine are exactly the same as in dimensions $24$ to $64$, with the same arguments. There are a few longer extensions:

\begin{prop} \label{etwentythreeextensions} $\;$

\begin{enumerate}
\item $e[124,6] \eta = e[125,21]$ \label{extensiononehundredandtwentyfour}
\item $e[129,3] \eta = e[130,18]$ \label{extensiononehundredandtwentynine}
\item $e[149,7] \eta = 2e[150,10] = e[150,22]$
\label{extensiononehundredandfortynine}
\end{enumerate}
\end{prop}

\begin{proof}
\begin{enumerate}
\item Due to the new $d_{23}$--differential and the $\nu$--extension \ref{efiveextensions}\eqref{firstnuextension}, we have that
\[
e[124,6] = \langle \copy\kbar^5,\copy\kbar,\nu \rangle,
\]
and multiplication with $\eta$ followed by Toda shuffling yields the extension to $e[125,21]$.
\item We have that $e[129,3] = \langle e[124,6],\eta,\nu\rangle$ and $e[130,18] = \langle e[125,21],\eta,\nu \rangle$ (the latter using \ref{efiveextensions}\eqref{epsilonextension}), thus this extensions follows from \eqref{extensiononehundredandtwentyfour}.
\item This follows from \eqref{extensiononehundredandtwentynine} by multiplying with $\copy \kbar$.\proved
\end{enumerate}
\end{proof}

There is one new differential, namely $d_{23}(e[146,2]) = e[145,25]$ which immediately follows from the $\eta$ extension \ref{etwentythreeextensions}\eqref{extensiononehundredandtwentyfour}.

\begin{figure}[ht!]
  \cl{\includegraphics{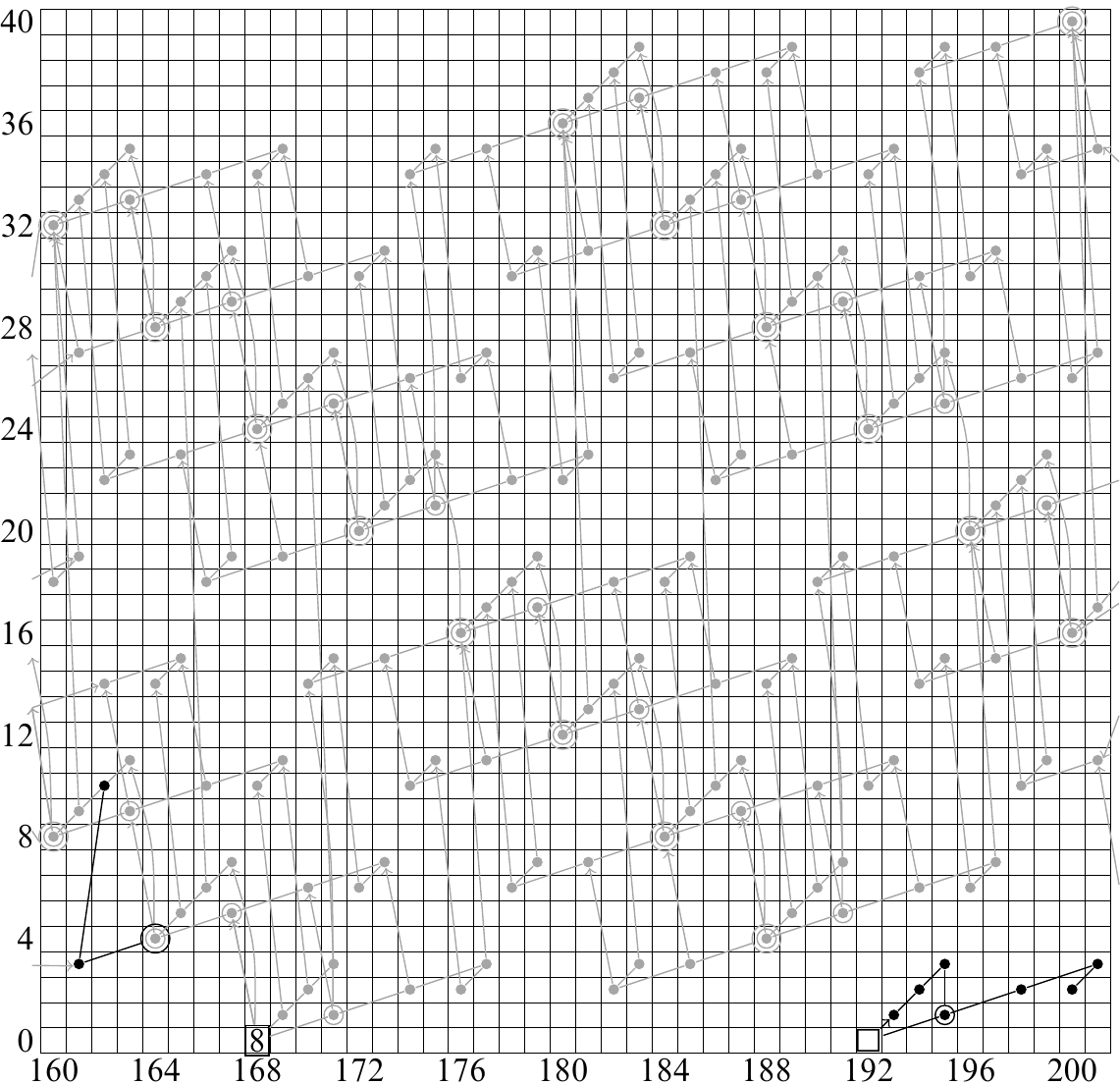}}
\caption{The elliptic spectral sequence --- dimensions 160--201}\label{chart160}
\end{figure}

The rest of the chart, up to the periodicity dimension of $192$, is displayed in \fullref{chart160}.

Since almost all classes die, there are only extensions that already occured in the respective dimensions lowered by $96 = |\Delta^4|$.

The new differential $d_{23}(e[171,3]) = e[172,26]$ follows from \ref{etwentythreeextensions}\eqref{extensiononehundredandfortynine}. From here on, all classes die until dimension $192$. This implies that $\Delta^8$ is a polynomial generator for the whole spectral sequence as well as in $\pi_*\tmf$.

\clearpage
\appendix
\section{Some Toda brackets in the homotopy of \texorpdfstring{$\tmf$}{tmf}} 
\setobjecttype{App}\label{todaapp}
In this appendix, I will assemble some $2$--primary Toda bracket
relations in $\tmf$ for quick reference. These are shown in \fullref{Table1}.
Many of them are proved somewhere in the text, others follow from the relation
\begin{equation}\label{todashuff}
\langle y,x,y \rangle \cap \langle x,y,2y \rangle \neq \emptyset
\end{equation}
valid in any $E_\infty$ ring spectrum for odd-dimensional classes $y$ 
(Toda \cite{Toda:composition}).

\begin{table}[h!]
\vglue0.5cm\small
\begin{tabular}{cccccc}
Dim & Class & $E_2$ repr. & Toda bracket & Massey product & Ref.\\
\hline
6 & $\nu^2$ & $h_2^2$ & $\langle \eta,\nu,\eta \rangle$ & $\langle h_1,h_2,h_1 \rangle$ & \eqref{masseyh22}\\
8 & $\epsilon$ & $c$ & $\langle \nu,\eta,\nu \rangle$ & $\langle h_2,h_1,h_2 \rangle$ & \eqref{masseyc}\\
&& & $\langle 2\nu,\nu,\eta \rangle$ & $\langle 2h_2,h_2,h_1 \rangle$ & \eqref{todashuff}\\
&&& $\langle \nu,2\nu,\eta \rangle$ & $\langle h_2, 2h_2,h_1 \rangle$ & juggling\\
&&& $\langle 2,\nu^2,\eta \rangle$ & $\langle 2,h_2^2,h_1 \rangle$ & juggling\\
14 & $\kappa$ & $d$ & $\langle \nu,2\nu,\nu,2\nu \rangle$ & $\langle h_2,2h_2,h_2,2h_2 \rangle$ & \eqref{masseyd}\\
15 & $\kappa \eta$ & $d h_1$ & $\langle \nu,\epsilon,\nu \rangle$ & $\langle h_2,c,h_2 \rangle$ & \eqref{masseydh1alt}, \eqref{todashuff}\\
&&& $\langle \epsilon,\nu,2\nu \rangle$ & $\langle c,h_2,2h_2 \rangle$ &\eqref{masseydh1alt}\\
20 & $\copy \kbar$ & $g$ & $\langle \kappa,2,\eta,\nu\rangle$ & --- & \eqref{kbartoda}\\
&&& $\langle \nu,2\nu,\nu,4\nu,\nu,\eta\rangle$ & --- & \fullref{kbarlemma}\\
&&& $\langle \nu,2\nu,\nu,4\nu,\eta,\nu\rangle$ & --- & \fullref{kbarlemma}\\
20 & $2 \copy \kbar$ & $2g$ & $\langle \kappa \eta,\eta,\nu \rangle$ & $\langle dh_1,h_1,h_2\rangle$ & \eqref{massey2g}\\
21 & $\copy \kbar \eta$ & $g h_1$ & $\langle \kappa,2,\nu^2 \rangle$ & --- & \eqref{kbaretatoda}\\
\\
\end{tabular}
\caption{Some Toda brackets and Massey products}\label{Table1}
\end{table}

We verify the remaining relations not explained in the main text:
\begin{align}
d h_1 &= h_1 \langle h_2,2h_2,h_2,2h_2 \rangle = \langle \langle h_1,h_2,2h_2 \rangle,h_2,2h_2 \rangle = \langle c,h_2,2h_2 \rangle \label{masseydh1alt}\\
\copy\kbar \eta &= \langle \kappa,2,\eta,\nu \rangle \eta = \langle \kappa,2,\langle \eta,\nu,\eta \rangle \rangle = \langle \kappa,2,\nu^2 \rangle \label{kbaretatoda}
\end{align}

\bibliographystyle{gtart}
\bibliography{link}

\end{document}